%% file: main.tex
\documentclass[a4paper,11pt]{article}

\usepackage{comment}
\usepackage{a4wide}
\usepackage[english]{babel}
\usepackage[utf8]{inputenc}
\usepackage[fleqn]{amsmath}
\usepackage{amsfonts}
\usepackage{amsthm}
\usepackage{bbm}
\usepackage{dsfont}
\usepackage{xcolor}
\usepackage{url}
\usepackage{enumitem}
\usepackage{hyperref}
\usepackage{mathtools}
\usepackage{tikz}
\usepackage{url}
\usepackage{nameref} 
\usepackage{hyperref} 
\usepackage{color}

\newtheorem{theorem}{Theorem}

\newtheorem{proposition}{Proposition}

\newtheorem{remark}{Remark}
\newtheorem{corollary}{Corollary}

\newcommand{\set}[1]{{\left\{ #1 \right\}}} % \set a b le pone { }
\newcommand{\inlinematrix}[4]{ \begin{pmatrix} #1 & #2 \\ #3 & #4 \end{pmatrix} }

\def\Z{\mathbb{Z}}

\allowdisplaybreaks

\title{Determinants of matrices related to the Pascal triangle}

\author{Martín Mereb}
\date{\today}

\def\modulo{\pmod}

\newcommand{\Ub}[4]{U_{#1,#2}(#3,#4)}

\newcommand{\U}[2]{U_{0,#2}(#1)}

\begin{document}

\maketitle

\begin{abstract}
In this note we prove an assertion made by 
M. Levin in 1999: the Pascal matrix modulo~$2$  has the property that each of the square sub-matrices laying on the upper border or on the  left border has determinants,   computed in~$\Z$, equal to~$1$ or~$-1.$ 
\end{abstract}

\section{Introduction}

 In this note we prove an that the Pascal matrix modulo~$2$  has the property that each of the square sub-matrices laying on the upper border or on the  left border has determinant,   computed in~$\Z$, equal to~$1$ or~$-1.$   
This extends some of the results in~\cite{Bacher2002, Bacher-Chapman} on determinants related to the Pascal triangle.

The Pascal triangle matrix 
%modulo $b\in\N$
has been  used in the theory of uniform distribution modulo one
%, first  by H. Faur\'e and   I.M. Sobol',  
to construct sequences of real numbers in the unit interval with  smallest possible discrepancy:
the first $N$ terms have discrepancy at most 
$(\log N) /N $ times a constant (see~\cite{Bugeaud2012} and the references therein).
When we restrict to sequences of  the form 
$\{b^n x \modulo 1\}_{n=1,2,3,\ldots}$ for any integer $b$ greater than or equal to~$2$ and for real numbers~$x$,  the smallest exact  discrepancy  that can be achieved by some~$x$  is not known. The question dates back to Korobov in 1956 (cfr.~\cite{Korobov1956}).
% https://en.wikipedia.org/wiki/Nikolay_Korobov

Using the Pascal triangle matrix modulo~$2$, 
in~\cite{Levin1999} M. Levin constructs   numbers 
$x$ such that the sequence 
$\{b^n x \modulo 1\}_{n=0,1,2,\ldots}$ has discrepancy of 
the first $N$ terms bounded from above by  $(\log N)^2/N $
times a constant. 
Becher and Carton in~\cite{BecherCarton2019}  defined variants of  the Pascal triangle matrix modulo~$2$  that  have the same property of the invertibility of the square sub-matrices  laying  on the upper or left border.  
They obtain a family of numbers  with the same property as Levin's. 
Larcher and Hofer recently showed that  for Levin's number constructed  for $b=2$ the discrepancy estimate $(\log N)^2/N $ is the best possible
(cfr.~\cite{HoferLarcher}).

The property that all square matrices in the upper and left border of the  Pascal matrix modulo 2  have determinants, computed in $\Z$ equal to $1$ or $-1$   ensures that if these determinants are computed in $\Z/b\Z$, for any $b$, they are also equal to $1$ or $-1$.
Thus, indeed, Levin's method yields 
numbers $x$ such that the 
$\{b^n x \modulo 1\}_{n=1,2,\ldots}$ has the 
small discrepancy property.

The article is organized as follows. In section~\ref{s:pascmat} we introduce some notation, define  the infinite matrix $U$ and state the main result (Theorem~\ref{t:main}). Section~\ref{s:proof} is devoted to its proof. In section~\ref{s:family} we define a whole family of matrices sharing the property of having all its sub-matrices laying on the upper or left border invertible, compute its number and give some examples.

%\section{Notation and Statement}
%\section{Preliminaries}
\section{Pascal matrices}\label{s:pascmat}

We want to study determinants of sub-matrices of certain infinite matrix. It will be convenient to index the rows and columns with non-negative numbers.

Let $U$ be the infinite matrix whose entry in 
the $i,j$ position is
the remainder when the binomial coefficient $\binom j i$ 
is divided by $2.$
%
%Here we use the standard notation for the binomial coefficient
%\begin{align*}    { \binom {n}{k} } := {\frac {n!}{k!(n-k)!} } \end{align*}
%when $ n \geq k,$ and $0$ otherwise.
%
%\begin{align*}
%U_{i,j} \equiv \binom{j}{i} \,\modulo 2.
%\end{align*}
%
Namely
\begin{align*}
U_{i,j} =
{\begin{cases}
1&{\text{ if }} \binom {j}{i} \text{ is odd, } \\
0&{\text{ if }} \binom {j}{i} \text{ is even.} 
\end{cases}}
\end{align*}

\begin{remark}\label{r:kummer}
By the well known result of Kummer~\cite{kummer} 
we know that the $(i,j)$-entry of $U$ is $1$
if $j \geq i$ and the binary representations of 
both $ i $ and $j-i$ don't share a $ 1 $ in the 
same position, and $0$ otherwise.

%\begin{align*}
%U[i][j] =
%{\begin{cases}
%1&{\text{if }}j \geq i \text{ and both },
%i \text{ and }j-i \text{ don't share a } 1 \text{ in the same %position in base }2 \\
%0&{\text{ otherwise.}}
%\end{cases}}
%\end{align*}
\end{remark}

% Notation for sub-matrices
Writing $ U _{I,J} $
for the sub-matrix of $U$ corresponding to the rows and columns indexed respectively by the sets 
$I,J \subseteq \Z_{\geq 0}$ introduce the following notation. 
The subset of integers $k$ greater than or equal to $m$
but smaller than $n$ is denoted by 
$[m:n) = \set{k \in \Z :, m\leq k < n}.$ 
The principal minors 
$U_{[0:n),[0:n)}$ 
%$\U n 0$
are denoted $U(n).$ 
For the top-most minors 
$U_{[0:n),[m:m+n)}$ 
%$\U n m$ 
we write  $\U n m .$ 
Finally, $\Ub i j n m $ stands for 
$U_{[i:i+n)[j:j+m)}.$
%$\Ub n m r k .$
%

The main result of this article is the following:
\begin{theorem}\label{t:main}
The sub-matrix 
%$U_{[0:n][m:m+n]}$ 
%$U(n,m)$ 
$\U n m $ 
has determinant $\pm 1,$
for all choices of $n,m\geq 0.$
\end{theorem}

\begin{remark}\label{r:tensor}
From the celebrated formula by Lucas~\cite{lucas}
it follows that $U$ is the infinite tensor of the matrix
$U(2)=%\smatrix 1 1 0 1 .$ $
\inlinematrix 1 1 0 1 .$
Here we think of this infinite tensor as the stable top-left square matrix in the sequence 
$U(2^k)= U(2)^{\otimes k }$ for $k\geq 1$. For example:
\begin{align*}    U(2) = &     \begin{pmatrix}     1 & 1 \\      0  & 1     \end{pmatrix} \\
%\end{align*}
%\begin{align*}
    U(4) %= U(2)^{\otimes 2} 
%    = U(2) \otimes U(2) 
    = 
    \begin{pmatrix} 
    U(2) & U(2) \\ 
    0  & U(2)
    \end{pmatrix}
    = &
    \begin{pmatrix} 
    1 & 1 & 1 & 1 \\ 
    0 & 1 & 0 & 1 \\
    0 & 0 & 1 & 1 \\
    0 & 0 & 0 & 1 
    \end{pmatrix}  
    \\
%\end{align*}
%\begin{align*}
    U(8) %= U(2)^{\otimes 3} 
%    = U(2) \otimes U(2)  \otimes U(2)  
    = 
%    \begin{pmatrix}    U(4) & U(4) \\     0  & U(4)     \end{pmatrix}     = 
    \begin{pmatrix} 
    U(2) & U(2) & U(2) & U(2) \\ 
    0 & U(2) & 0 & U(2) \\
    0 & 0 & U(2) & U(2) \\
    0 & 0 & 0 & U(2) 
    \end{pmatrix} 
    =&
    \begin{pmatrix} 
    1 & 1 & 1 & 1 & 1 & 1 & 1 & 1 \\  
    0 & 1 & 0 & 1 & 0 & 1 & 0 & 1 \\
    0 & 0 & 1 & 1 & 0 & 0 & 1 & 1 \\
    0 & 0 & 0 & 1 & 0 & 0 & 0 & 1 \\
    0 & 0 & 0 & 0 & 1 & 1 & 1 & 1 \\  
    0 & 0 & 0 & 0 & 0 & 1 & 0 & 1 \\
    0 & 0 & 0 & 0 & 0 & 0 & 1 & 1 \\
    0 & 0 & 0 & 0 & 0 & 0 & 0 & 1 
    \end{pmatrix} 
    \\
\vdots &
\end{align*}

\end{remark}

\begin{remark}\label{r:periodic}
The first $2^k$ rows of $U$ are $2^k$-periodic
%(n < 2^k)
.
\end{remark}

\begin{remark}\label{r:uptriang}
All the sub-matrices 
$U(n)%=U_{[0:n][0:n]}
$ 
are upper triangular with only 1's on the diagonal.
\end{remark}

\begin{remark}\label{r:symmetry}
$U(2^k) % = U(2)^{\otimes k }
%U_{[0:2^{k}][0:2^{k}]}
$
are symmetric with respect to the anti-diagonal.
\end{remark}

%\section{Proof of Main result}
\section{Proof of Theorem~\ref{t:main} }\label{s:proof}

%\begin{proof}[]
\proof[Proof of Theorem~\ref{t:main}]
We proceed by induction on the size $n$ of the sub-matrix.
Since the first row of $U$ is made entirely of $1'$s,
the statement is true for $n = 0.$
By Remark~\ref{r:uptriang} the statement is also 
true for $m = 0.$
Consider $k$ such that $2^{k-1} < n \leq 2^k.$
By periodicity (cfr. Remark~\ref{r:periodic})
we may assume $m < 2^k.$
We separate in two cases according to whether
$m$ is less than $2^{k-1}$ or not.
%\begin{enumerate}[label=(\roman*)]
%\item 
\subsection*{{Case $m < 2^{k-1}$:}}

%\begin{flushleft}\begin{large}{\textbf{{Case $m < 2^{k-1}$:}}}\end{large}\end{flushleft}

Let us compare the matrices
%$U_{[0:n][m:m+n]}$  
%$U(n,m)$  
$\U n m$ 
and  
%$U_{[0:n][m+2^{k-1}:m+2^{k-1}+n]}.$ 
%$U(n,m+2^{k-1}).$ 
$\U n {m+2^{k-1}}$ 
The first $2^{k-1}$ rows are identical by 
$2^{k-1}$-periodicity (cfr. Remark~\ref{r:periodic}).

For the remaining $n-2^{k-1}$ rows we apply 
elementary row operations 
to the first matrix and obtain the second one,
up to sign.

Subtracting
%$U_{[0:n-2^{k-1}),[m:m+n)}$  
$\Ub 0 m {n-2^{k-1}} n $  
from
%$U_{[2^{k-1}:n),[m:m+n)}  $  
$\Ub {2^{k-1}} m {n-2^{k-1}} n $  
%one obtains
we get exactly 
(by Remark~\ref{r:kummer})
the sub-matrix 
%$U_{[2^{k-1}:n),[m+2^{k-1}:m+2^{k-1}+n)}$ 
$\Ub {2^{k-1}} {m+2^{k-1}} {n-2^{k-1}} n $  
multiplied by $ -1$ (Fig.~\ref{fig:casesmall}).

\input{casesmall}

Therefore, the determinant of 
%$U(n,m)$ 
$\U n m $ 
%$U_{[0:n)[m:m+n)}$  
is 
that of 
%$$U_{[0:n)[m+2^{k-1}:m+2^{k-1}+n)}$$ 
%$U(n,m+2^{k-1})$ 
$\U n {m+2^{k-1}}$ 
multiplied by 
$(-1)^{(n-2^{k-1})}$
and this case reduces to the next one.
%\item 
\subsection*{{Case $m \geq 2^{k-1} $:}}

%\begin{flushleft}\begin{large}{\textbf{Case $m \geq 2^{k-1} $:}}\end{large}\end{flushleft}

Subdivide 
$\U n m$ 
%$U(n,m)$ 
%$U_{[0:n)[m:m+n)}$ 
into $4$ blocks as follows (Fig.~\ref{fig:caselarge}):
%
%\begin{align*}     U(n,m) = \begin{pmatrix}    U_{[0:m+n-2^{k})[m:2^{k})} &    U_{[0:m+n-2^{k})[2^{k}:m+n)} \\ U_{[m+n-2^{k}:n)[m:2^{k})} & U_{[m+n-2^{k}:n)[2^{k}:m+n)}    \end{pmatrix}. \end{align*}
%
taking 
%$I=[0:m+n-2^{k})$, \\ $I'=[m+n-2^{k}:n) $,  
%$J=[m:2^{k})$ and $J'=[2^{k}:m+n)$ 
$J=[m:2^{k})$, $J'=[2^{k}:m+n)$, 
$I=[0:m+n-2^{k})$ and $I'=[m+n-2^{k}:n) $,  
we get the partition
\begin{align*}
%    U_{[0:n)[m:m+n)} = 
%    U( n, m) = 
    \U n m = 
%    U_{I\cup I', J\cup J'} = 
    \begin{pmatrix}
    U_{I,J}  &
    U_{I,J'}  \\ 
    U_{I',J}  & 
    U_{I',J'}
    \end{pmatrix}
\end{align*}
as ${[0:n)} = I\cup I'$ and ${[m:m+n)} = J\cup J'.$

Note that the 
bottom-right block $U_{I',J'}$ is full of zeros 
(by $2^k-$periodicity and Remark~\ref{r:uptriang})
and the 
top-right one $U_{I,J'}$ agrees with
$\U {m+n-2^{k}} {2^{k}} $
%$U(m+n-2^{k}, 2^{k})$
%$U_{[0:m+n-2^{k})[2^{k}:m+n)}$
(again by Remark~\ref{r:periodic})
and therefore has determinant $1.$

To find the determinant of 
%$U(n,m)$ 
$\U n m$ 
we consider the
block matrix obtained by
swapping the blocks
\begin{align*}
    \begin{pmatrix}
    U_{I,J'}  &
    U_{I,J}  \\ 
    U_{I',J'}  & 
    U_{I',J}
    \end{pmatrix}
%=    \begin{pmatrix}    U_{[0:m+n-2^{k})[2^{k}:m+n)} &    U_{[0:m+n-2^{k})[m:2^{k}) }\\     U_{[m+n-2^{k}:n)[2^{k}:m+n)} &    U_{[m+n-2^{k}:n)[m:2^{k})}     \end{pmatrix}.
\end{align*}
which is upper triangular by blocks
with first block $ U_{I,J'}$ having determinant $1.$

\input{caselarge}

The determinant of 
$\U n m$ 
%$ U(n,m) $
%$U_{[0:n)[m:m+n)}$ 
ends up being 
%times 
$(-1)^{(n-1)(m+n-2^{k} )} =
(-1)^{(n-1)m} $
times
that of 
%$U_{[m+n-2^{k}:n),[m:2^{k})}$
$U_{I',J}$.

Note that by Remark~\ref{r:symmetry} 
the \emph{anti-transpose} 
$U_{I,J}^\ast$
of this last sub-matrix\footnote{i.e.: the reflection with respect to the anti-diagonal.}
$U_{I',J}$ is 
%$U_{[0:2^{k}-m)[2^{k}-n: 2^{k} - n + 2^{k}-m )}.$
%$U_{I',J}^\ast = U(2^{k}-m, 2^{k}-n)$ 
\\
$U_{I',J}^\ast = \U {2^{k}-m}  {2^{k}-n} $ 
(Fig.~\ref{fig:anti}).

\input{antitranspose}

Since
$$
2^{k}-m \leq 2^{k}-2^{k-1} = 2^{k-1} < n
$$
its determinant is $\pm 1$ by inductive hypothesis.
\endproof

%\end{enumerate}

%\end{proof}

%

%\section{Additional comments}

\section{A family of Pascal-like matrices}\label{s:family}

%\begin{remark} 
Our matrix $U$ appeared in~\cite{Bacher-Chapman}
to prove that all symmetric Pascal matrices $\modulo{2}$ have determinant $\pm 1.$

Consider $P$ the infinite 
matrix\footnote{it is noted $\overline{P}(\infty)_2$ in~\cite{Bacher-Chapman}.}
\begin{align*}
P_{i,j} =
{\begin{cases}
1&{\text{ if }} \binom {j+i}{i} \text{ is odd, } \\
0&{\text{ if }} \binom {j+i}{i} \text{ is even.} 
\end{cases}}
\end{align*}
For its $LU$-decomposition we define $L$ as the transpose of $U$ and $D$ as the infinite diagonal with 
the Thue--Morse sequence
\footnote{\href{https://oeis.org/A106400}{https://oeis.org/A106400}}. Namely
\begin{align*}
    D_{i,j} =
{\begin{cases}
1&{\text{ if }} i=j \text{ has an even number of $1$'s in base $2$, } \\
-1&{\text{ if }} i=j \text{ has an odd number of $1$'s in base $2$, } \\
0&{\text{ if }} i\neq j. 
\end{cases}}
\end{align*}

Adopting the same notation for sub-matrices as with $U$
we have 
$P , L$ and $D$ are the infinite tensor product of
\begin{align*}
P(2) = \inlinematrix 1 1 1 0 , \qquad 
L(2) = \inlinematrix 1 0 1 1 \qquad \text{ and } \qquad 
D(2) = \inlinematrix 1 0 0 {-1} 
%, \qquad U(2) = \inlinematrix 1 1 0 1
\end{align*}
respectively.

Since 
$$\inlinematrix 1 1 1 0  = \inlinematrix 1 0 1 1
\inlinematrix 1 0 0 {-1} \inlinematrix 1 1 0 1$$
we get that $P(n) = L(n)D(n)U(n)$ for all $n\geq 1.$

As $L$ is lower triangular and $D$ is diagonal with only $\pm 1$'s, Theorem~\ref{t:main} %immediately 
implies
\begin{corollary}\label{c:main}
Given $n,m\geq 0,$ the sub-matrix 
%$P_{[0:n][m:m+n]}$ 
%$P(n,m)$ 
$P_m (n)$
has determinant $\pm 1.$
\end{corollary}

The same result applies for infinite matrices having a similar $LU$-decomposition.
%\end{remark}

By the symmetry of $P$ we deduce that every square sub-matrix laying on the upper or left border of $P$
has determinant $\pm 1.$ 

We say that a matrix is \emph{Pascal-like} if
every such sub-matrix is invertible. When working over the integers this means having determinant $\pm 1$ so
our matrix $P$ is Pascal-like by Corollary\ref{c:main}.

Another example over the integers is provided by the honest Pascal matrix $M_{i,j} = \binom{i+j}{j}$ as can be seen by a routine application of Vandermonde determinant and elementary row operations.

More is true according to the following
\begin{proposition}\label{p:pascal-like}
Let $R$ be a commutative ring with finite group of units $R^\times.$
The number of Pascal-like matrices $M\in R^{n\times m}$ is exactly $\#(R^\times)^{nm}.$
\end{proposition}
\proof
Each entry $M_{i,j}$ is the bottom-right entry of exactly one square sub-matrix laying on the top or left side whose determinant is to be a unit.
There are precisely $\#(R^\times)^{nm}$ ways to prescribe those determinants. For each such prescription there is a unique way of solving for each
entry $M_{i,j}$ recursively in $i+j$ by row expansion.
\endproof

\begin{corollary}
There are exactly $2^{nm}$ Pascal-like matrices in $\Z^{n\times m}$. All of them are congruent $\modulo 2$
but at the same time they cover all the possible Pascal-like matrices when reduced modulo $3,4$ or $6.$
\end{corollary}
\proof
This follows from Proposition~\ref{p:pascal-like} and the facts that
$$\#(\Z^\times) = \varphi(3) = \varphi(4) = \varphi(6) = 2$$ 
and $\varphi(2) = 1.$
\endproof

\section*{Acknowledgements}

The author would like to thank Prof. Becher for bringing this problem to his attention in the first place, together with a lot of helpful comments as well.

\bibliographystyle{plain}
\bibliography{refs.bib}

\noindent
Martín Mereb \\
 Departamento de Matemática, Facultad de Ciencias Exactas y Naturales \& IMAS \\
 Universidad de Buenos Aires \&  CONICET Argentina-  {\tt  mmereb@gmail.com}
\end{document}

%% file: casesmall.tex
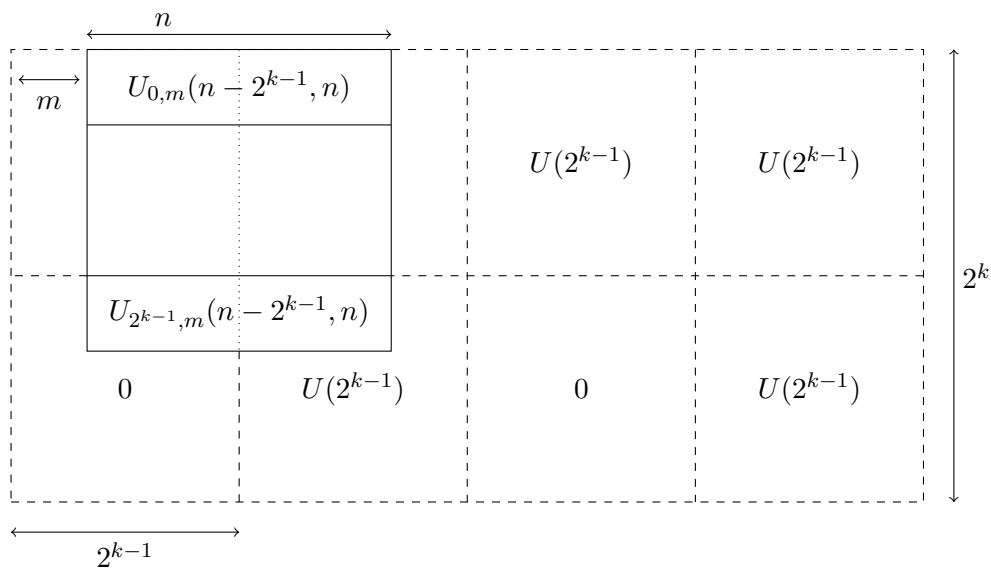
\begin{figure}[htbp]
\begin{center}
    \begin{tikzpicture}%[transform canvas={scale=.5}]
        \draw [dashed](-6,3) rectangle (6,-3);
        \draw [dashed] (0,3) -- (0,-3);
        \draw [dashed] (3,3) -- (3,-3);
        %\draw [dashed] (-3,3) -- (-3,-3);
        \draw [dotted] (-3,3) -- (-3,-1);
        \draw [dashed] (-3,-1) -- (-3,-3);
        \draw [dashed] (-5,0) -- (-6,0);
        \draw [dashed] (6,0) -- (-1,0);
        \draw (-5,0) -- (-1,0);
        \draw (-5,2) -- (-1,2);
%        \node at (-3,2.5) {$U_{[0:n-2^{k-1}),[m:m+n)}$};
        \node at (-3,2.5) {$\Ub 0 m {n-2^{k-1}} n $};
        \node at (-3,-0.5) 
        {$\Ub {2^{k-1}}{m}{n-2^{k-1}}{n}$};
%        \node at (-3,-0.5) {$U_{[2^{k-1}:n),[m:m+n)}$};
        
        \draw [<->] (6.4,3) -- (6.4,-3);
        \node at (6.7,0) {$2^k$};
        
        \draw [<->] (-6,-3.4) -- (-3,-3.4);
        \node at (-4.5,-3.7) {$2^{k-1}$};
        
        \draw (-5,3) rectangle (-1,-1);
        \draw [<->] (-5,3.2) -- (-1,3.2);
        \node at (-4,3.4) {$n$};
        \draw [<->] (-5.9,2.6) -- (-5.1,2.6);
        \node at (-5.5,2.3) {$m$};
        
        %\node at (-4,1.5) {$U_{I,J}$};
        %\node at (-4,-0.5) {$U_{I',J}$};
        %\node at (-2,1.5) {$U_{I,J'}$};
        %\node at (-2,-0.5) {$U_{I',J'}$};
        
        \node at (1.5,1.5) {$U(2^{k-1})$};
        \node at (4.5,1.5) {$U(2^{k-1})$};
        \node at (4.5,-1.5) {$U(2^{k-1})$};
        \node at (1.5,-1.5) {$0$};
        
        \node at (-4.5,-1.5) {$0$};
        \node at (-1.5,-1.5) {$U(2^{k-1})$};

    \end{tikzpicture}
    \end{center}
        \caption{Position of 
%        $U(n,m)$ 
        $\U n m $ 
        when $m < 2^{k-1}$.}
    \label{fig:casesmall}   
\end{figure}

%% file: caselarge.tex
\begin{figure}[htbp]
\begin{center}
    \begin{tikzpicture}%[transform canvas={scale=.5}]
        \draw [dashed](-6,3) rectangle (6,-3);
        \draw [dashed] (0,-1) -- (0,-3);
        \draw [dashed] (-3,3) -- (-3,-3);
        \draw [dashed] (3,3) -- (3,-3);
        \draw [dashed] (6,0) -- (2,0);
        \draw [dashed] (-2,0) -- (-6,0);
        \draw [dotted] (-2,0) -- (2,0);
        
        \draw [<->] (6.4,3) -- (6.4,-3);
        \node at (6.7,0) {$2^k$};
        
        \draw [<->] (-6,-3.4) -- (-3,-3.4);
        \node at (-4.5,-3.7) {$2^{k-1}$};
        
        \draw (-2,3) rectangle (2,-1);
        \draw [<->] (-2,3.2) -- (2,3.2);
        \node at (0,3.4) {$n$};
        \draw [<->] (-5.9,2.6) -- (-2.1,2.6);
        \node at (-2.5,2.3) {$m$};
        
        \draw (-2,1) -- (2,1);
        \draw (0,3) -- (0,-1);
        
        \draw [dotted] (0,3) -- (6,-3);
        
        \draw [<->] (0,-3.2) -- (2,-3.2);
        \node at (1,-3.5) {$m+n-2^{k}$};

        \node at (-1,2) {$U_{I,J}$};
        \node at (-1,0) {$U_{I',J}$};
        \node at (1,2) {$U_{I,J'}$};
        \node at (1,0) {$U_{I',J'}$};
        
        \node at (-1.5,-1.5) {$U(2^{k-1})$};
        \node at (1.5,-1.5) {$0$};
        \node at (-4.5,1.5) {$U(2^{k-1})$};
        \node at (4.5,1.5) {$U(2^{k-1})$};
        \node at (4.5,-1.5) {$U(2^{k-1})$};
        \node at (-4.5,-1.5) {$0$};
        
    \end{tikzpicture}
    \end{center}
        \caption{Position of 
%        $U(n,m)$ 
        $\U n m$ 
        when $m\geq 2^{k-1}$.}
    \label{fig:caselarge}   
\end{figure}
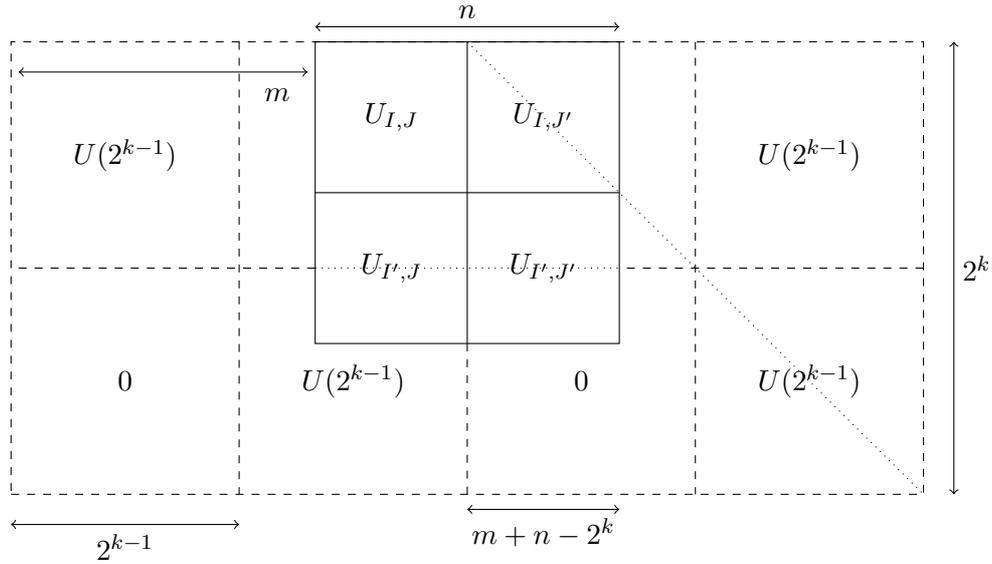

%% file: antitranspose.tex
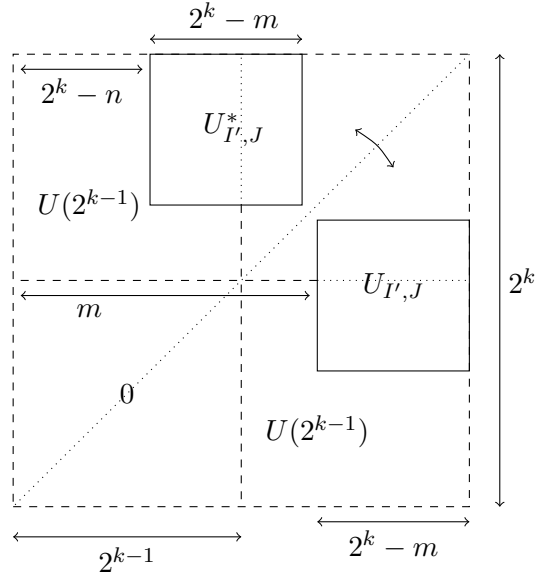
\begin{figure}[htbp]
\begin{center}
    \begin{tikzpicture}%[transform canvas={scale=.5}]
        \draw [dashed](-3,3) rectangle (3,-3);
        \draw [dashed] (0,1) -- (0,-3);
        \draw [dotted] (0,1) -- (0,3);
        \draw [dashed] (1,0) -- (-3,0);
        \draw [dotted] (1,0) -- (3,0);

        \draw [dotted] (3,3) -- (-3,-3);
        
        \draw [<->] (3.4,3) -- (3.4,-3);
        \node at (3.7,0) {$2^k$};
        \draw [<->] (-3,-3.4) -- (0,-3.4);
        \node at (-1.5,-3.7) {$2^{k-1}$};

        \draw (1,0.8) rectangle (3,-1.2);
        
        \draw [<->] (-2.9,2.8) -- (-1.3, 2.8);
        \node at (-2.1,2.5) {$2^k - n$};
        
        \draw [<->] (-1.2,3.2) -- (0.8,3.2);
        \node at (-0.1,3.5) {$2^{k}-m$};
        
        \draw (0.8,1) rectangle (-1.2,3);
        
        \draw [<->] (-2.9,-.2) -- (0.9,-.2);
        \node at (-2,-.4) {$m$};
        
        \draw [<->] (1,-3.2) -- (3,-3.2);
        \node at (2,-3.5) {$2^{k}-m$};
        
        \node at (2,-0.1) {$U_{I',J}$};
        
        \node at (-0.1,2) {$U_{I',J}^\ast$};
        
        %\node at (2,2) {$U(2^{k-1})$};
        \draw[<->] (2,1.5) .. controls (1.8,1.8) and (1.8,1.8) .. (1.5,2);
        \node at (-2,1) {$U(2^{k-1})$};
        \node at (1,-2) {$U(2^{k-1})$};
        \node at (-1.5,-1.5) {$0$};
        
    \end{tikzpicture}
    \end{center}
    \begin{center}
        \caption{
        %The sub-matrix 
%        $U(2^k-m,2^k-n)$ coincides with
        $\U {2^k-m} {2^k-n}$ coincides with
        $U_{I',J}^\ast,$ \\ the anti-transpose of $U_{I',J}$.}        
    \end{center}
    \label{fig:anti}   
\end{figure}

%% file: main.bbl
\begin{thebibliography}{1}

\bibitem{Bacher2002}
Roland Bacher.
\newblock Determinants of matrices related to the {P}ascal triangle.
\newblock {\em J. Th\'{e}or. Nombres Bordeaux}, 14(1):19--41, 2002.

\bibitem{Bacher-Chapman}
Roland Bacher and Robin Chapman.
\newblock Symmetric {P}ascal matrices modulo {$p$}.
\newblock {\em European J. Combin.}, 25(4):459--473, 2004.

\bibitem{BecherCarton2019}
Ver\'{o}nica Becher and Olivier Carton.
\newblock Normal numbers and nested perfect necklaces.
\newblock {\em J. Complexity}, 54:101403, 12, 2019.

\bibitem{Bugeaud2012}
Yann Bugeaud.
\newblock {\em Distribution modulo one and {D}iophantine approximation}, volume
  193 of {\em Cambridge Tracts in Mathematics}.
\newblock Cambridge University Press, Cambridge, 2012.

\bibitem{HoferLarcher}
Roswitha Hofer and Gerhard Larcher.
\newblock The exact order of discrepancy for levin's normal number in base 2.
\newblock May 2022.

\bibitem{Korobov1956}
Nikolai~Mikhailovich Korobov.
\newblock On completely uniform distribution and conjunctly normal numbers.
\newblock {\em Izv. Akad. Nauk SSSR. Ser. Mat.}, 20:649--660, 1956.

\bibitem{kummer}
Ernst~Eduard Kummer.
\newblock \"{U}ber die {E}rg\"{a}nzungss\"{a}tze zu den allgemeinen
  {R}eciprocit\"{a}tsgesetzen.
\newblock {\em J. Reine Angew. Math.}, 44:93--146, 1852.

\bibitem{Levin1999}
Mordekhaĭ~Borisovich Levin.
\newblock On the discrepancy estimate of normal numbers.
\newblock {\em Acta Arith.}, 88(2):99--111, 1999.

\bibitem{lucas}
Édouard Lucas.
\newblock Sur les congruences des nombres eul\'{e}riens et les coefficients
  diff\'{e}rentiels des functions trigonom\'{e}triques suivant un module
  premier.
\newblock {\em Bull. Soc. Math. France}, 6:49--54, 1878.

\end{thebibliography}
